\documentclass{ifacconf}

\usepackage{graphicx}      % include this line if your document contains figures
\usepackage{natbib}        % required for bibliography
\usepackage{amsfonts, amsmath, amssymb}
\usepackage{tikz}
\usepackage{tikz-network}
\usepackage{subcaption}
\usepackage{todonotes}
\usepackage{enumerate}
\usepackage{float}
\usepackage{xcolor}

\newtheorem{definition}{Definition}

\newtheorem{theorem}{Theorem}

\newtheorem{example}{Example}

\newtheorem{assumption}{Assumption}

\AtEndEnvironment{example}{\hfill$\lozenge$}
\AtEndEnvironment{theorem}{\hfill$\triangle$}

\global\long\def\G{\mathcal{G}}
\global\long\def\E{\mathcal{E}}
\global\long\def\V{\mathcal{V}}

\global\long\def\R{\mathbb{R}}

\global\long\def\ker{\mathrm{Ker}}

\graphicspath{{figures/}}

\begin{document}
\begin{frontmatter}

  \title{A geometric view of formation control with application to
  directed sensing\thanksref{footnoteinfo}}

  \thanks[footnoteinfo]{This work was supported by the Israel Science
    Foundation grant no. 453/24 and UK Research and Innovation though
    the grant UKRI1112.
    We thank the ICMS in Edinburgh, where this project was initiated,
  for its hospitality. }

  \author[LST]{Louis Theran}
  \author[DZ]{Daniel Zelazo}
  \author[JS]{Jessica Sidman}

  \address[LST]{School of Mathematics and Statistics,
  University of St Andrews, Scotland (e-mail: lst6@st-andrews.ac.uk).}
  \address[DZ]{Stephen B. Klein Faculty of Aerospace Engineering,
    Technion - Israel Institute of Technology, Haifa, Israel (e-mail:
  dzelazo@technion.ac.il)}
  \address[JS]{Department of Mathematics, Amherst College, Amherst, MA, USA
  (e-mail: jsidman@amherst.edu)}

  \begin{abstract}                % Abstract of 50--100 words
    We propose a geometric approach to distance-based formation control
    modeled on a minimum-norm lifting of
    Riemannian gradient descent in edge-space to node-space.  This
    yields a unified family of controllers, including the classical
    gradient controller and its directed variant.
    For the directed case, we give a simple numerical test
    for local convergence that applies to any directed graph and target.
    We show that  persistence is neither necessary nor sufficient for
    local convergence of our directed controller and propose an
    alternative that is necessary and more easily checked.
  \end{abstract}


\end{frontmatter}
%===============================================================================

\section{Introduction}

The distance-constrained formation control problem has emerged as one
of the canonical problems for multi-agent coordination.  The
objective is to steer a group of agents to a prescribed geometric
shape, specified by inter-agent distances, using only local sensing
and distributed control actions.  The classical solution,
by \cite{krickBrouckeFrancis}, is a
gradient-based controller that uses the squared edge length
error as a potential in node space and follows the negative
gradient.  When the sensing graph is generically $d$-rigid,
the gradient controller is generically locally exponentially
stable, as shown by \cite{Sun2016}, using tools from
rigidity theory.  In a notable work, \cite{Dorfler_TAC2010}
introduced a differential-geometric perspective to the
problem that indicates the squared edge lengths evolve on a
manifold of feasible distances in edge space.

A more realistic model of the problem is the
\emph{directed} variant, in which each distance constraint
is assigned to one agent.  The resulting interactions become
non-reciprocal and lose their potential-gradient structure.  As a
result, directed formation control is significantly
less well-understood.  Much of the existing literature
focuses on specific graph families derived from
\emph{persistence}, which is a property of a directed graph that
attempts to adapt rigidity theory to the directed
setting motivated by the role of $d$-rigidity in
the solution to the undirected problem (\cite{krickBrouckeFrancis}
and Theorem \ref{thm: model}); see also
\cite{Hendrickx_IJRNC2007,Anderson_CDC2007,Yu_SIAMJCO2009,Babazadeh_CDC2020}.
A complementary line of work by \cite{Zhang_JDSMC2018}
introduces gain design to recover stability
with directed sensing; they are able to do this for the
highly restricted class of persistent graphs in
leader / first-follower form.

In this work, we take a different perspective on the directed
problem.  Instead of working with a property of oriented graphs,
we work with $d$-rigid graphs and impose a condition on the quadratic
form derived from a modification of
the rigidity matrix that incorporates the data of the graph
orientation (Theorem~\ref{thm: dircontrol}).
This connects the information from the directed graph to the
manifold of feasible distance measurements, a core object of
study in rigidity theory.
It also allows us to
identify open, semialgebraic regions of the target space on
which the dynamics are local exponentially stable, and to
provide combinatorial certificates when no such region exists
(Theorem \ref{thm: admiss}).
Summarizing, we initiate a new, structural theory of robust
directed formation control.

The approach to the directed formation control problem described
above rests on  taking a differential-geometric viewpoint of the
distance-based formation control focused on the
dynamics of an artificial edge-space system evolving on
the manifold of feasible measurements.  Riemmannian
gradient descent on the edge manifold, followed by a minimum
norm lift to node dynamics gives a geometrically canonical
``model'' solution to formation control
(Section \ref{sec.formationProblem}). While not practically
useful in a distributed setting, the model controller
serves as a template for a family of  controllers
that includes the model, gradient, and directed
controllers (Section \ref{sec.mainresult}).  We derive a
simple condition, that can be checked with linear algebra,
for local exponential stability in this
family (Theorem \ref{thm: main}) that implies generic
exponential convergence of the model and gradient controllers
(Theorem \ref{thm: model}).  Then we specialize to the
directed controller, giving a simple, numerical condition sufficient
for local exponential convergence that applies to any
graph and target formation (Theorem \ref{thm: dircontrol}).
Simulations in Section \ref{sec.examples}
validate our theory and illustrate that local exponential stability is
not a generic property of the directed controller, showing that
persistence is neither necessary nor sufficient in the directed
setting.  This leads us to introduce the notions of dynamic and algebraic
admissibility, which are generic properties of a directed graph that can
be easily checked, and which we prove to be
necessary for local exponential stability
of the directed controller (Theorem \ref{thm: admiss})
in Section \ref{sec.combinatorics}.

\paragraph*{Notations:}
We write \(\mathbb{R}\) for the set of real numbers and
\(\mathbb{Q}\) for the set of rational numbers. For a matrix \(A\),
we write \(A^\top\) for its transpose and
\(A^\dagger\) for its Moore--Penrose pseudoinverse. We use \(A \succ
0\) (\(A \succeq 0\)) to indicate that \(A\) is symmetric positive
definite (semidefinite).
The image, kernel, and rank of \(A\) are written \(\mathrm{Im}\,A\),
\(\ker\, A\), and \(\mathrm{rk}\,A\).
For a subspace \(W\), its orthogonal complement is \(W^\perp\).
For a smooth map \(F\), the differential at \(x\) is denoted \(\mathrm{d}F_x\).
All inner products are Euclidean unless stated otherwise, and
\(\|\cdot\|\) denotes the associated norm for vectors, and induced
norm for operators. For a manifold \(\mathcal
M\), \(T_x\mathcal M\) denotes its tangent space at \(x\) and
\(\mathrm{grad}_{\mathcal M} V\) denotes the Riemannian gradient of a
smooth function \(V\colon \mathcal M \to \mathbb{R}\).  The quadratic
form of an operator $T$ is denoted $q_T(x) := \langle x,T(x)\rangle$,
and the restriction of $T$ to a subspace $V$ is denoted $T|_V$.

\section{The Formation Control Problem}\label{sec.formationProblem}

In this section, we briefly review the standard formulation of the
formation control problem and then introduce a closely related,
artificial problem that exposes its underlying geometric structure.

The starting point for this work is with an ensemble of $n$ agents,
modeled by integrator dynamics,
\begin{align}\label{integrator}
  \dot p_i(t) & = u_i(t),\quad i=1,\ldots,n,
\end{align}
where $p_i(t) \in \R^{d}$ is the state (position) of agent $i$, and
$u_i(t) \in \R^d$ is the control.  The aggregate position (control)
vector is denoted as $p(t) =
\begin{bmatrix} p_1(t)^\top & \cdots & p_n(t)^\top
\end{bmatrix}^\top$ (similarly for $u(t)$). We assume the agents
interact by an information exchange network described by the
undirected graph~$\G=(\V,\E)$.

The \emph{target configuration} $p^*$ has an associated set
of desired inter-agent square distances, $m_{ij}^*$, for each $ij
\in \E$; $m^*
\in \R^{|\E|}$ is the aggregate target distance vector.
We  define the \emph{distance map} $F \,:\,
\mathbb{R}^{d n} \to \mathbb{R}^{ |\E|}$ by
\begin{align}\label{distance_map}
  F(p) =
  \begin{bmatrix} \cdots & \|p_j-p_i\|^2 & \cdots
  \end{bmatrix}^\top,
\end{align}
so that $F(p^*) = m^*$.  From now on, the squared edge length
measurement of a target formation $p^*$ will be denoted $m^*$.
The distance-based formation control
problem is to drive an initial configuration $p^0$ to a point $q$
congruent and close to $p^*$ (such that $F(q)=m^*$),
using, at each time $t$, the information $m^*$ and  $m(t) = F(p(t))$.
We introduce two
assumptions that we will use at various points in this paper.
\begin{assumption}\label{asump.drigid-reg}
  The graph $\mathcal G$ is generically $d$-rigid
  and the target configuration $p^*$ is a regular point of the rigidity map $F$.
\end{assumption}
\begin{assumption}\label{asump.drigid-geneic}
  The graph $\mathcal G$ is generically $d$-rigid
  and the target configuration $p^*$ is generic.
\end{assumption}
Here generic means that $p^*$ holds over the
open, dense complement of an unspecified exceptional set
that is defined by polynomials with rational coefficients.
Every generic configuration is regular, and if $p^*$ is
generic, then $F(p^*)$ is generic in the image of $F$.
A graph $\mathcal{G}$ is $d$-rigid if, for every configuration
$p^*$ so that $F(p^*)$ is smooth in the image of $F$, the
framework $(\mathcal{G},p^*)$ is locally rigid; in particular, when
$\mathcal{G}$ is $d$-rigid, the frameworks $(\mathcal{G},p^*)$
for regular $p^*$ will be rigid (\cite{Gortler_AJM2010}).

\subsection{The Standard Control Law}
In the standard formulation of the formation control problem, a
common approach is to use a
negative-gradient flow of the potential function $V(p) =
\tfrac{1}{4}\|F(p)-m^*\|^2$, yielding the control
\begin{align}\label{grad_Control}
  u(t) &= -\tfrac{\partial V(p(t))}{\partial p(t)}=
  -R(p(t))^\top(R(p(t))p(t)-m^*).
\end{align}
Here, we note that $R(p(t))p(t) = F(p(t))$, where $R(p(t))\in
\mathbb{R}^{n d \times |\mathcal E|}$ is known as the \emph{rigidity matrix}.
This control ensures local exponential convergence to a configuration $q$
near $p^*$ such that $F(q) = m^*$,
which, if $\mathcal{G}$ is $d$-rigid will imply that $q$ is
congruent to $p^*$ (see, e.g., \cite{krickBrouckeFrancis, Sun2016,
Dorfler_TAC2010} and Theorem \ref{thm: model}).

\subsection{The Edge System}

The work of \cite{Dorfler_TAC2010} analyzed \eqref{grad_Control} from
a geometric perspective.  The starting point of that analysis was to
consider the dynamics of distances, i.e., $\tfrac{d}{dt}F(p(t))$.  In
this work we take a similar approach, although our starting point
will not be the agent level dynamics defined by \eqref{grad_Control}.
Rather, we will assume the existence of an artificial system of edges
with state $m(t) \in \R^{|\E|}$ that evolve according to the integrator dynamics
\begin{align}\label{edge_integrator}
  \dot m(t) = v(t),
\end{align}
for some control $v(t)$.  In the ``edge space" we now seek a control
to drive the edge states to the desired set of square distances $m^*$.

We emphasize that, unless $n\le d + 1$, there must be algbraic
constraints on $v(t)$ to ensure that $m(t)$ evolves in a
way that is consistent with a physical system; i.e., that there is a
$p(t)$ such that $F(p(t)) = m(t)$ for all $t$.
Since the constraint is that $m(t)$ lies in the image of $F$, we
denote by $\mathcal{Q}$ the regular values of
the map $F$.  Because $F$ is a polynomial mapping defined over
$\mathbb{Q}$, the semialgebraic Sard theorem \cite[p. 234]{RAG}
implies that $\mathcal{Q}$ is a smooth
semialgebraic set that is open and dense in the image of $F$.
(The image of $F$ itself is not smooth.  We will see presently why
we use $\mathcal{Q}$ instead of the smooth locus of ${\rm Im}\, F$.)
We also note that $\mathcal{C} :=
F^{-1}(\mathcal{Q})\subseteq \R^{dn}$, the set of
regular points of $F$, is a smooth semi-algebraic set of dimension $dn$.
At any point $m = F(p)$ in $\mathcal{Q}$,
the tangent space $T_m\mathcal Q$ is the vector space of differential
changes to the squared
edge lengths that can arise from the infinitesimal motions of
the agents, i.e.,
\[
  T_m\mathcal{Q} = \mathrm{Im}\,R(p).
\]
Accordingly, the edge dynamics~\eqref{edge_integrator} are to be
interpreted as evolving on $\mathcal{Q}$,
meaning that the velocity $v(t)$ must satisfy $v(t)\in
T_{m(t)}\mathcal{Q}$ for all $t$.

The manifold $\mathcal Q$ inherits the Euclidean inner product from
$\mathbb{R}^{|\mathcal E|}$, defining a Riemannian metric,
\[
  g_m(\xi_1,\xi_2) = \xi_1^\top \xi_2, \qquad \xi_1,\xi_2\in T_m\mathcal Q.
\]
Let $\Pi(m)\,:\,\R^{|\mathcal{E}|}\to
T_m\mathcal Q$ denote the orthogonal projector onto $T_m\mathcal Q$.
If $m$ is a regular value of $F$, for any $p\in F^{-1}(m)$, we have
\begin{equation}\label{projector}
  \Pi(m) = R(p)R^\dagger(p).
\end{equation}
(Note that if $m$ is smooth but not regular,
  then the formula for $\Pi(m)$ will be valid for some
$p$ in the fiber over $m$ but not others.)

Define the \emph{edge potential} function,
\[
  V_e(m) = \tfrac12\|m - m^*\|^2,
\]
which measures the squared deviation of the current
edge vector $m$ from the desired one $m^*$.
As we are interested in feasible motions that must remain
on $\mathcal Q$, the appropriate notion of the gradient is the
\emph{Riemannian gradient} on $(\mathcal Q,g)$, obtained by
projecting $\nabla V_e(m)= m - m^*$ onto the tangent space of $\mathcal{Q}$
at $m$,
\[
  \nabla_{\!\mathcal Q}V_e(m) = \Pi(m)\,\nabla V(m)
  = \Pi(m)\,(m - m^*).
\]

The \emph{Riemannian gradient flow} associated with $V_e$ is defined by
\begin{equation}\label{eq:riem-flow}
  \dot m = -\,\nabla_{\!\mathcal Q}V_e(m)
  = -\,\Pi(m)\,(m - m^*).
\end{equation}
A solution to \eqref{eq:riem-flow} describes a trajectory $m(t)\in\mathcal
Q$ whose velocity $\dot m(t)$ is the projection of the negative
Euclidean gradient of $V_e$ onto the tangent space $T_{m(t)}\mathcal Q$.
Hence, the motion remains on $\mathcal Q$ for all time and follows
the steepest descent of $V_e$ with respect to the Riemannian metric~$g$.

\subsection{Connecting Edge Flows to Node Flows}\label{subsec.modelcontrol}
To connect the squared edge length flow defined by the Riemmanian
gradient flow in \eqref{eq:riem-flow} to node dynamics, we make use of the
fact that $\mathcal{Q} = F(\mathcal{C})$.
At any configuration $p\in \mathcal{C}$, the differential of the edge
map $F$ is the
map $${\rm d}F_p : T_p\mathcal{C}\longrightarrow T_{F(p)}\mathcal
Q,$$ where we recall that ${\rm d}F_p =  2R(p)$.
This linear map relates infinitesimal changes of the node positions $\dot p=u$
to infinitesimal changes of the squared edge lengths via
\[
  \dot m = {\rm d}F_p\, u = 2R(p)\,u.
\]
Hence, to produce a prescribed edge-space velocity $v^*\in
T_{F(p)}\mathcal Q$, we seek
a node-space velocity $u$ satisfying $${\rm d}F_p\,u=v^*.$$

Among the infinitely many $u$ that satisfy this relation,
a natural choice is the one with \emph{minimum Euclidean norm}:
\begin{align}\label{eq:node-opt}
  \min_{u\in\mathbb{R}^{dn}} \;\|u\|^2
  \quad\text{s.t.}\quad {\rm d}F_p\,u = v^*.
\end{align}
The unique minimizer of~\eqref{eq:node-opt} is obtained via the
Moore--Penrose pseudoinverse,
\[
  u^* = ({\rm d}F_p)^\dagger v^* = 2\,R(p)^\dagger v^*.
\]
This solution corresponds to the \emph{minimum-norm configuration
velocity} that generates
the desired edge-space motion $v^*$.

The mapping ${\rm d}F_p^\dagger$ provides the least-energy lift of an
admissible edge velocity
to a consistent motion of the agent configuration.
Consequently, the combined edge--node system can be viewed as a
two-level geometric control
construction:
\begin{enumerate}
  \item In edge space, compute the Riemannian gradient flow
    $v^*=-\nabla_{\!\mathcal Q}V_e(m)$,
    which defines the optimal tangent velocity on~$\mathcal Q$.
  \item In node space, realize this motion through the minimum-norm lift
    $u^* = {\rm d}F_p^\dagger v^*$.
\end{enumerate}
Substituting $v^*=-\,\Pi(m)\,(m-m^*)$ yields the node-space dynamics
\begin{align}\label{model_controller}
  \dot p = u^*
  = -\,\tfrac12\,R(p)^\dagger\,\Pi(m)\,(m-m^*).
\end{align}
This expression represents the \emph{minimum-energy configuration flow}
that is consistent with the Riemannian steepest-descent dynamics on
the edge manifold.

We refer to \eqref{model_controller} together with
\eqref{eq:riem-flow} as the \emph{model controller}.  This control,
while not distributed as in \eqref{grad_Control}, provides a clear
geometric picture of what should be done to drive agents to a desired formation.

\section{From the Model Controller to Distributed
Solutions}\label{sec.mainresult}
While the model controller \eqref{eq:riem-flow} and
\eqref{model_controller} represents a geometrically
canonical approach to solving the formation control problem,
it is not intended as an algorithm for practical
implementation.  Instead, we use it as a template for
a new design space of controllers with the same essential
characteristics:
edge dynamics evolving on $\mathcal{Q}$ and
lifted node dynamics.  The larger class is big
enough to include distributed architectures such as
\eqref{grad_Control} and controlled enough to be
amenable to combinatorial analysis.

\subsection{Edge-driven controller family}
Now we introduce the control family we study in this section.
The state space for the controllers is
\[
  \mathcal{S} = \{(p,m)\in \mathcal{C}\times \mathcal{Q} : F(p) = m\},
\]
which is a smooth, irreducible, semi-algebraic set equipped with the natural
projections $\pi_1$ and $\pi_2$.
At any point $(p,m) \in \mathcal S$, its tangent space is
\[
  T_{(p,m)} \mathcal{S} = \{(\dot p,\dot m) : \dot m = 2 R(p)\dot p\},
\]
which is the same coupling of node and edge velocities that we have seen in
the derivation of the model law \eqref{model_controller} and
\eqref{eq:riem-flow}.  We will study control laws determined
by state-dependent linear transformations
\[
  \nu_{(p,m)} : \R^{|\mathcal{E}|}\to (\R^d)^n \qquad \text{and}\qquad
  \eta_{(p,m)} : \R^{|\mathcal{E}|}\to \R^{|\mathcal{E}|},
\]
such that $(\nu_{(p,m)} x,\eta_{(p,m)} x)\in T_{(p,m)} \mathcal{S}$
for all states
in $\mathcal{S}$ and vectors $x\in \R^{|\mathcal{E}|}$.
Given a target
formation $(\mathcal{G},p^*)$ with edge measurements
$m^*$ and a starting formation $(\mathcal{G},p^0)$, the
controller is the initial value problem (IVP)
\begin{subequations}\label{eq:pm-dynamics}
  \begin{align}
    \dot p &= \nu_{(p,m)}(m^* - m) \label{eq:pm-dynamics-a}\\
    \dot m &= \eta_{(p,m)}(m^* - m) \label{eq:pm-dynamics-b}
  \end{align}
\end{subequations}
with initial conditions
\begin{align*}
  \begin{cases}
    \dot p(0) &= \nu_{(p^0,F(p^0))}(m^* - F(p^0))\\
    \dot m(0) &=  \eta_{(p^0,F(p^0))}(m^* - F(p^0)).
  \end{cases}
\end{align*}
For technical reasons, we will also allow controllers to be defined only
on an open neighborhood in $\mathcal{S}$ of $(p^*,m^*)$.  Controllers
in this family yield solutions $(p(t),m(t))\in \mathcal{S}$.  When
we discuss the \emph{edge dynamics} of one of these controllers, we
mean the flow $\pi_2(p(t),m(t))\in \mathcal{Q}$, and the node dynamics are,
similarly, $\pi_1(p(t),m(t))\in \mathcal{C}$.  The edge and node
dynamics have the
same smoothness as the full solution.
For notational ease, we use the shorthand $\eta := \eta_{(p,m)}$
and $\eta^*:=\eta_{(p^*,m^*)}$ (similarly for $\nu$). In the
following examples we illustrate how the model and gradient
controllers fit into this setup.  With an eye to applying
algebraic geometry, we
show that the transformations $\eta_{(p,m)}$ and $\nu_{(p,m)}$ are
rational functions of the state.

\begin{example}[Model controller, revisited]\label{exm: model law}
  The model controller \eqref{eq:riem-flow} and
  \eqref{model_controller} arises from setting
  \[
    \nu = \frac{1}{2}R^\dagger(p)\quad \text{and}\quad
    \eta = \Pi(m)= 2R(p)\nu.
  \]
  The formulas for $\eta$ and $\nu$ contain the
  pseudo-inverse of $R(p)$, which is analytic in the state, but not
  rational.  However, in a neighborhood of $p^*$, we can define
  the projectors as rational functions if we fix the indices of
  the columns of
  a basis for the column space of $R(p)$.
  Thus, there is a neighborhood of $(p^*,m^*)$
  on which the model law is defined by rational functions.
\end{example}

\begin{example}[Gradient controller, revisited]\label{exm: grad law}
  The classical gradient law arises from setting
  \[
    \nu = R^\top(p)\text{ and } \eta = 2R(p)R^\top(p).
  \]
  These are  polynomial functions of the state $(p,m)$, which are
  rational over all of  $\mathcal{S}$. The column space of
  $\eta$ is contained in the column space of $R(p)$, and,
  hence, in $T_m \mathcal{Q}$.  By construction, we have
  $2R(p)\nu = \eta$,
  so the controller is well-defined.
\end{example}
We note that the matrix $\eta$ in Example \ref{exm: grad law} is
not a projection, since it does not fix its image pointwise.

\subsection{Local exponential stability}
We now provide a sufficient condition for the local exponential
stability of \eqref{eq:pm-dynamics} to a
generic target.
\begin{theorem}[Local exponential stability]\label{thm: main}
  Let Assumption \ref{asump.drigid-reg} hold.
  Suppose
  that for some neighborhood $V\subseteq \mathcal{S}$ of $(p^*,
  m^*)$, the linear transformations $\nu$ and $\eta$ are rational
  functions of $(p,m) \in V.$  If the quadratic form $q_{\eta^*}$
  restricted to ${T_{m^*} \mathcal{Q}}$,
  is positive definite,
  then there is a neighborhood $U\ni (p^*,m^*)$ such that the control
  law determined by
  \eqref{eq:pm-dynamics-a}--\eqref{eq:pm-dynamics-b} with initial condition
  $(p^0,F(p^0)) \in U$ converges to a state $(q,m^*)$ with
  $q$ congruent to $p^*$.  Moreover,
  the edge dynamics on $\mathcal{Q}$ and the node dynamics on
  $\mathcal{C}$ are locally exponentially stable.
\end{theorem}
We remark that, in the proof of Theorem \ref{thm: main}, given in
Sec.~\ref{subsec.thm1proof}, the assumption that the controller is
rational can be
substantially relaxed, to, e.g., $C^2$.  With rationality,
we get, additionally, the algebro-geometric notion of genericity,
which corresponds to structural robustness of the control law.
The local exponential convergence certificate given
by the positivity of $q_{\eta^*}$ will, if it holds at one
point, hold on an open semialgebraic subset of target configurations.
Hence, non-robust stability is \emph{non-generic}, and the
sufficient condition
provided in Theorem \ref{thm: main} is robust to perturbation in the
target geometry. For the model and gradient controllers, the situation
is even better, and local convergence is a generic property.
Although the sufficient stability condition in
Theorem \ref{thm: main} is strictly
stronger than the linearized dynamics being Hurwitz, for these
two controllers,
it is also necessary, yielding an extremely clean, combinatorial theory.
\begin{theorem}[Generic undirected exponential stability]\label{thm: model}
  Let Assumption \ref{asump.drigid-reg} hold.  Then there exists
  a neighborhood $U$ of
  $p^*$ such that the node dynamics on $\mathcal{C}$
  of the model controller \eqref{model_controller},
  \eqref{eq:riem-flow}, and the gradient controller
  \eqref{grad_Control}  onto $\mathcal C$
  converge exponentially  to a configuration congruent to $p^*$
  for any initial
  condition $p\in U$  if and only if $\mathcal{G}$ is $d$-rigid.
\end{theorem}
\begin{pf}
  The case where $\mathcal{G}$ is not $d$-rigid is easier, so we
  do it first.
  If $\mathcal{G}$ is not $d$-rigid, because $p^*$
  generic, every neighborhood of $p^*$ contains a configuration $q$
  not congruent to $p^*$ such that $F(p^*) = F(q)$.  At this
  $q$, both the model and gradient controllers will immediately stop
  (i.e., $\dot p = 0$).  Hence,
  they do not locally converge to a configuration congruent to $p^*$.

  Now suppose that $\mathcal{G}$ is $d$-rigid.
  The semialgebraic Sard Theorem \cite[p. 234]{RAG} implies that
  $p^*$ is a regular point
  of $F$, so $(p^*,m^*)\in \mathcal{S}$, putting us in the general setting
  of Theorem \ref{thm: main}, which we use to show local exponential
  convergence.
  Examples  \ref{exm: model law} and \ref{exm: grad law}
  show that the model and gradient controllers are rationally defined and
  satisfy the other hypotheses of Theorem \ref{thm: main}.  What remains is
  to check that $q_{\eta^*}$ is positive definite.  In both laws,
  the $\eta_{(p,m)}$
  matrix is PSD with columns spanning $T_m \mathcal{Q}$.  Since the
  restriction of a
  PSD matrix to its column space is positive definite, the result follows.
  The verification of PSD in the model controller is because $\Pi(m)$ is a
  projector and in the gradient controller because $\eta = 2R(p)R(p)^\top$
  is presented as factoring.    \hfill $\qed$
\end{pf}
The conclusion for the gradient controller is due to
\cite{krickBrouckeFrancis}
and has been reproved many times in the literature.

\subsection{Proof of Theorem \ref{thm: main}}\label{subsec.thm1proof}
Let us fix a target $(p^*,m^*)\in \mathcal{S}$.  The goal is to
find a neighborhood $U\subseteq \mathcal{S}$ of $(p^*,m^*)$ so that
there are:
\begin{enumerate}[(I)]
  \item Constants $C,c > 0$ such that, for all initial conditions
    $(p^0,m^0)\in U$, and all $t\ge 0$,
    $\|m^* - m(t)\| \le C{\rm e}^{-ct}$;
  \item Constants $K,k > 0$ such that, for all initial conditions
    $(p^0,m^0)\in U$, there is  a configuration $q$
    congruent to $p^*$, such that
    $\|q - p(t)\| \le K{\rm e}^{-kt}$.
\end{enumerate}
The proof is geometric and illustrates many of the ideas leading to
the derivation of the model controller and larger family.  Along the
way, we mention connections to standard techniques in control theory
and rigidity.

\subsubsection{Node stability from edge stability}
We first assume that (I) holds and derive (II).  The essence of
the argument is that smoothness of $\nu$ implies that the edge-to-node
dynamics are input-to-state stable (ISS).  Rigidity enters here
because, locally,
the fiber $F^{-1}(m^*)$ consists of configurations congruent to
$p^*$.  Now we give the formal argument.

Assume (I), i.e., there are constants $C,c > 0$ and a neighborhood
$V\subseteq \mathcal{Q}$
of $m^*$ so that, for all initial conditions $(p^0,m^0)\in
F^{-1}(V)\times V$,
$\|m^* - m(t)\| \le C{\rm e}^{-ct}$.

Because $\mathcal{G}$ is
$d$-rigid and $p^*$ is a regular value of $F$,
we can find a neighborhood $U_0$ of $p^*$ so that
if $q\in U$ and $F(q) = m^*$, then $q$ is congruent
to $p$.  Let us now consider the full dynamics in the
neighborhood $U = (U_0\cap F^{-1}(V)) \times V$ of the
target $(p^*,m^*)$.  Shrinking $U$, if necessary,
we can find an $N > 0$ so that $\|\nu\| < N$,
because $\nu$ is rational in the state,
and so its eigenvalues  remain  bounded on a neighborhood of $(p^*,m^*)$.
Using our
estimate on $\|\nu\|$, we obtain
\[
  \|\dot{p}(t)\|= \|\nu(m^*-m)\|\le N\|m^*-m\| \le
  \underbrace{NC}_{C'}e^{-ct}.
\]
Consider a sequence $t_i \to \infty$, and let $j > i$.  We have
\[
  \|p(t_j) - p(t_i)\| \le \int_{t_i}^{t_j} \|\dot p(t)\|dt
  \le \int_{t_i}^\infty  C'e^{-ct} \le \frac{C'}{c}e^{-ct_i}.
\]
Hence, the sequence $p(t_i)$ is Cauchy, and has a limit $q$.
A similar estimate give us
\[
  \|q - p(t)\| \le \frac{C'}{c}e^{-ct},
\]
so the convergence is exponential.  Since the whole state
converges to $(q,m^*)\in \mathcal{S}$, we have $F(q) = m^*$ and
$q\in U_0$, so
$q$ is congruent to $p^*$. This completes the proof of (II)
assuming (I).  We remark that, had we not assumed exponential convergence,
some additional effort would be required to show the node dynamics do not
become trapped in the space of rigid motions.

\subsubsection{Exponential stability in edge space}
This is the core technical argument, and it proceeds parallel to
a standard line in Lyapunov stability theory.  The key novelty is that,
using the geometric perspective,
we identify a geometrically canonical quadratic energy,
namely $q_{\eta^*}$ on $T_{m^*}\mathcal{Q}$,
that certifies local exponential convergence in a way that can be
verified computationallyﬁ.

The proof of (I) will follow once we have a neighborhood
$V\subseteq \mathcal{Q}$ of $m^*$ so that, for $m\in V$
we have the coercive inequality
\begin{equation}\label{eq: coerce}
  q_{\eta}(m^* - m)\ge \alpha \|m^* - m\|^2,
\end{equation}
for some $\alpha > 0$.
Using
\eqref{eq: coerce}, we get
\[
  \frac{d}{dt} \|m^* - m(t)\|^2 =-2q_{\eta}(m^*-m) \le -2\alpha
  \|m^* - m\|^2,
\]
from which a standard application of Grönwall's inequality
\cite[Lemma A.1]{khalil} yields
exponential convergence of any solution of the edge dynamics with
initial condition in $V$.  The
same estimate shows that $(p,m)\mapsto (\nu,\eta)$
is a smooth vector field on $F^{-1}(V)\times V$, so the solution
exists, and the statement for edge convergence holds.  Hence,
we have (I), up to the construction of $V$.

By compactness and the hypothesis that $\eta^*$ induces a positive
definite quadratic form
on $T_{m^*} \mathcal{Q}$, there is a $\beta > 0$ so that
$q_{\eta^*}(x) > \beta$ for all
unit vectors $x$ in $T_{m^*} \mathcal{Q}$.

For each $m\neq m^*$, we write $e = \frac{m^* - m}{\|m^*-m\|}$ and then
$e = v + w$, where $v\in T_{m^*}\mathcal{Q}$ and $w\in
(T_{m^*}\mathcal{Q})^{\perp}$.
Because $\mathcal{Q}$ is an embedded smooth manifold,
for all $\varepsilon > 0$ there is an $r > 0$
so that, when $\|m^* - m\| < r$, $\|v\| \ge 1 - \varepsilon$ and
$\|w\| < \varepsilon$.  After repeated applications of this fact and
smoothness of
$\eta$ in the state, we find a neighborhood $V$ of $m^*$ on which
for all $m\in V$ and $(p,m)\in F^{-1}(V)\times V$,
\begin{align}
  \|v\| > & 1/2 \label{eq: est1}\\
  \|w\| < & \beta/(32\|\eta^*\|) \label{eq: est2} \\
  \|w\|^2 < &\beta/(16\|\eta^*\|) \label{eq: est3} \\
  \|\eta^* - \eta\| < & \beta/16 \label{eq: est4}.
\end{align}
This is possible because the constants on the right are
independent of $m$ and
smoothness of $\eta$ in the state.  We claim this $V$ works.
We first lower bound, for $m\in V$,
\[
  q_{\eta^*}(m^* - m) = \|m^* - m\|^2 q_{\eta^*}(e).
\]
From the decomposition of $e$, we have
\[
  q_{\eta^*}(e) = q_{\eta^*}(v) + w^\top\eta^*v+v^\top\eta^*w +
  q_{\eta^*}(w).
\]
From the \eqref{eq: est1},
\[
  q_{\eta^*}(v) > \beta/4.
\]
From \eqref{eq: est2} and  $\|v\|\le \|e\| = 1$, we get
\[
  |w^\top\eta^*v|+|v^\top\eta^*w| \le 2\|v\|\|w\|\|\eta^*\| < \beta/16
\]
and, using \eqref{eq: est3}
\[
  q_{\eta^*}(w) \le \|w\|^2\|\eta^*\| < \beta/16.
\]
Combining these bounds yields
\[
  q_{\eta^*}(e) \ge \beta/4 -2\beta/16 =  \beta/8.
\]
Hence, for $m\in V$,
\[
  q_{\eta^*}(m^* - m)  \ge (\beta/8)\|m^* - m\|^2.
\]
Finally, by \eqref{eq: est4},
\[
  q_\eta(e) \ge (\beta/8 - \|\eta^* - \eta\|)\|e\|^2 \ge
  \beta/16\|m^* - m\|^2,
\]
which proves \eqref{eq: coerce} with $\alpha = \beta/16$.  This completes
the proof of (II).
\hfill $\qed$

\section{Directed Control Laws}\label{sec.directed}
An appealing feature of the gradient control \eqref{grad_Control} is
that it admits a \emph{distributed} implementation.  At the agent level,
\eqref{grad_Control} is equivalent to
\begin{equation}\label{agent_grad_control}
  u_i = \sum_{ij\in \E} (\|p_i-p_j\|^2-m^*)(p_j-p_i).
\end{equation}
Each agent therefore only requires information from its neighboring
agents defined by the
graph $\G$.  The model controller \eqref{model_controller},
however, leads to a
\emph{centralized} architecture, since the pseudo-inverse of the
rigidity matrix will generally be dense.

While distributed architectures are advantageous,
even the control \eqref{agent_grad_control} has shortcomings.
Indeed, the control assumes bidirectional information exchange.  In
many real-world applications, it is more realistic (especially with
on-board sensing) that the information exchange should be
\emph{directed}, i.e., agent $i$ sensing agent $j$ does not imply
agent $j$ is
sensing agent $i$.  We model this asymmetry by orienting the edges of the
sensing graph $\mathcal{G}$.
Given an orientation $\overrightarrow{\mathcal{G}}$ on the edges of
$\mathcal{G}$, we may simply  modify
\eqref{agent_grad_control}, as in \cite{Zhang_JDSMC2018}, to
\begin{equation}\label{dir_agent_grad_control}
  u_i = \sum_{\overrightarrow{ij}\in \E} (\|p_i-p_j\|^2-m^*)(p_j-p_i),
\end{equation}
where $\overrightarrow{ij}$ is the directed edge from node $i$ to
$j$.  This leads to the
aggregate dynamics, which is no longer a
gradient flow, of
\begin{equation}\label{dir_formation_ctr}
  \dot p = -\overrightarrow R(p(t))^\top(R(p(t))p(t)-m^*),
\end{equation}
where $\overrightarrow{R}(p)$ is obtained from $R(p)$ by zeroing out, in the
row corresponding to $ij$, the entries
corresponding to the node $j$ for the orientation $\overrightarrow{ij}$
and the entries corresponding to $i$ for the orientation
$\overleftarrow{ij}$.
We now describe a numerical test that is sufficient for
\eqref{dir_formation_ctr}
to converge to a given configuration $(\G, p^*)$.  Since it is based on
Theorem \ref{thm: main}, this result inherits the same structural stability
properties.
\begin{theorem}[Exponential stability of
  \eqref{dir_formation_ctr}]\label{thm: dircontrol}
  Let Assumption \ref{asump.drigid-reg} hold, and $(p^*,m^*)\in
  \mathcal{S}$ be the
  target state.  Let ${\small \overrightarrow{\mathcal{G}}}$ be
  an orientation of $\mathcal{G}$.
  Let $P$ be a matrix with
  linearly independent columns spanning the column space of $R(p^*)$
  and $Z = 2R(p)\,\overrightarrow{R}(p)^\top$.
  If $\frac{1}{2}P^\top (Z + Z^\top)P$ is positive definite, then
  then trajectories of \eqref{dir_formation_ctr} starting
  sufficiently close to
  $p^*$ converge exponentially to a configuration $q$ with $F(q)=m^*$.
\end{theorem}
\begin{pf}
  The edge dynamics for the directed law \eqref{dir_formation_ctr}
  are given by
  \[
    \dot m
    = 2R(p)\left[-\overrightarrow{R}(p)^\top(m-m^*)\right].
  \]
  So, in the setting of Theorem \ref{thm: main}, the directed control
  law corresponds to setting
  \[
    \nu := \overrightarrow{R}(p)^\top,
    \text{ and }
    \eta := 2R(p)\,\overrightarrow{R}(p)^\top = Z.
  \]
  The stated convergence result is an application of Theorem
  \ref{thm: main}, so
  we check the hypotheses.  That $\mathcal{G}$ is generically
  $d$-rigid has been
  assumed, as has $(p^*,m^*)\in \mathcal{S}$.
  The symmetric matrix $\frac{1}{2}P^\top (Z + Z^\top)P$ corresponds to the
  quadratic form of $\eta^*$ restricted to $T_{m^*}\mathcal{Q}$,
  so, when it is positive definite Theorem \ref{thm: main} applies, and
  we conclude the convergence result.
  \hfill $\qed$
\end{pf}
Even though the directed controller does not have a gradient-potential
structure, the stability certificate in Theorem \ref{thm: dircontrol} is a
natural geometric counterpart to the gradient potential for the directed
setting.

\section{Numerical Studies}\label{sec.examples}
We illustrate the behavior of the model controller
\eqref{model_controller}, the gradient controller
\eqref{grad_Control}, and the directed controller
\eqref{dir_formation_ctr} on representative examples.
\begin{figure}
  \begin{minipage}[b]{0.15\textwidth}
    \centering
    \begin{tikzpicture}
% Axis lines
    \draw[->, thin] (-1.3,0) -- (1.5,0) node[below right] {}; % x-axis
    \draw[->, thin] (0,-1.3) -- (0,1.5) node[above] {}; % y-axis
    
    % central node
    \Vertex[size=.3, x=0, y=0, label=$1$, style={shape=circle}]{v1}

    % outer nodes (4-cycle)
    \Vertex[size=.3, x=-.5,  y=-.5,  label=$2$, style={shape=circle}]{v2}
    \Vertex[size=.3, x=-1,  y=1,  label=$3$, style={shape=circle}]{v3}
    \Vertex[size=.3, x=2/3, y=1,  label=$4$, style={shape=circle}]{v4}
    \Vertex[size=.3, x=1,  y=-1, label=$5$, style={shape=circle}]{v5}

    % rim edges (cycle on 2–3–4–5)
    \Edge (v2)(v3)
    \Edge (v3)(v4)
    \Edge (v4)(v5)
    \Edge (v5)(v2)

    % spokes (center to rim)
    \Edge (v1)(v2)
    \Edge (v1)(v3)
    \Edge (v1)(v4)
    \Edge (v1)(v5)

\end{tikzpicture}
    % \\[0.3em]
    (a) $(\mathrm W_5,p^*)$ %Undirected wheel graph $\mathrm W_5$.
  \end{minipage}
  %\hfill
  \begin{minipage}[b]{0.15\textwidth}
    \centering
    \begin{tikzpicture}

% Axis lines
    \draw[->, thin] (-1.3,0) -- (1.5,0) node[below right] {}; % x-axis
    \draw[->, thin] (0,-1.3) -- (0,1.5) node[above] {}; % y-axis
    
    % central node
    \Vertex[size=.3, x=0, y=0, label=$1$, style={shape=circle}]{v1}

    % outer nodes (4-cycle)
    \Vertex[size=.3, x=-.5,  y=-.5,  label=$2$, style={shape=circle}]{v2}
    \Vertex[size=.3, x=-1,  y=1,  label=$3$, style={shape=circle}]{v3}
    \Vertex[size=.3, x=2/3, y=1,  label=$4$, style={shape=circle}]{v4}
    \Vertex[size=.3, x=1,  y=-1, label=$5$, style={shape=circle}]{v5}

    % rim edges (cycle on 2–3–4–5)
    \Edge[Direct](v2)(v3)
    \Edge[Direct](v3)(v4)
    \Edge[Direct](v4)(v5)
    \Edge[Direct](v5)(v2)

    % spokes (center to rim)
    \Edge[Direct](v1)(v2)
    \Edge[Direct](v1)(v3)
    \Edge[Direct](v4)(v1)
    \Edge[Direct](v5)(v1)
\end{tikzpicture}
    % \\[0.3em]
    (b) $(\overrightarrow{\mathrm W}_5, p^*)$ %Directed wheel graph
    % %$\overrightarrow{\mathrm W}_5$.
  \end{minipage}
  \begin{minipage}[b]{0.15\textwidth}
    \centering
    \begin{tikzpicture}

% Axis lines
    \draw[->, thin] (-1.3,0) -- (1.5,0) node[below right] {}; % x-axis
    \draw[->, thin] (0,-1.3) -- (0,1.5) node[above] {}; % y-axis
    
    % central node
    \Vertex[size=.3, x=9/5, y=-5/3, label=$1$, style={shape=circle}]{v1}

    % outer nodes (4-cycle)
    \Vertex[size=.3, x=-.5,  y=-.5,  label=$2$, style={shape=circle}]{v2}
    \Vertex[size=.3, x=-1,  y=1,  label=$3$, style={shape=circle}]{v3}
    \Vertex[size=.3, x=2/3, y=1,  label=$4$, style={shape=circle}]{v4}
    \Vertex[size=.3, x=1,  y=-1, label=$5$, style={shape=circle}]{v5}

    % rim edges (cycle on 2–3–4–5)
    \Edge[Direct](v2)(v3)
    \Edge[Direct](v3)(v4)
    \Edge[Direct,bend=-15](v4)(v5)
    \Edge[Direct,bend=-15](v5)(v2)

    % spokes (center to rim)
    \Edge[Direct,bend=15](v1)(v2)
    \Edge[Direct,bend=-15](v1)(v3)
    \Edge[Direct,bend=15](v4)(v1)
    \Edge[Direct](v5)(v1)
\end{tikzpicture}
    % \\[0.3em]
    (c) $(\overrightarrow{\mathrm W}_5, q^*)$%Directed wheel graph
    % $\overrightarrow{ \mathrm W}_5$.
  \end{minipage}
  \caption{(a) and (b)
    have the same target configuration; (b) and (c) have the same
    orientation,
    but $p^*$ satisfies Theorem~\ref{thm: dircontrol} while $q^*$ does not.
  }
  \label{fig:wheel5_combined}
\end{figure}
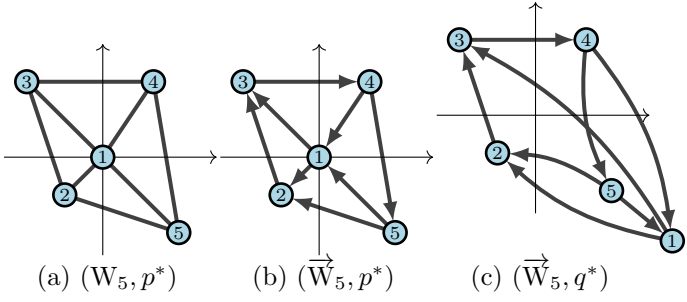
\subsubsection*{Comparison of the model and gradient controllers}
\begin{figure}[H]
  \centering
  % (a) Edge trajectories
  \begin{subfigure}[t]{0.24\textwidth}
    \centering
    \includegraphics[width=\linewidth]{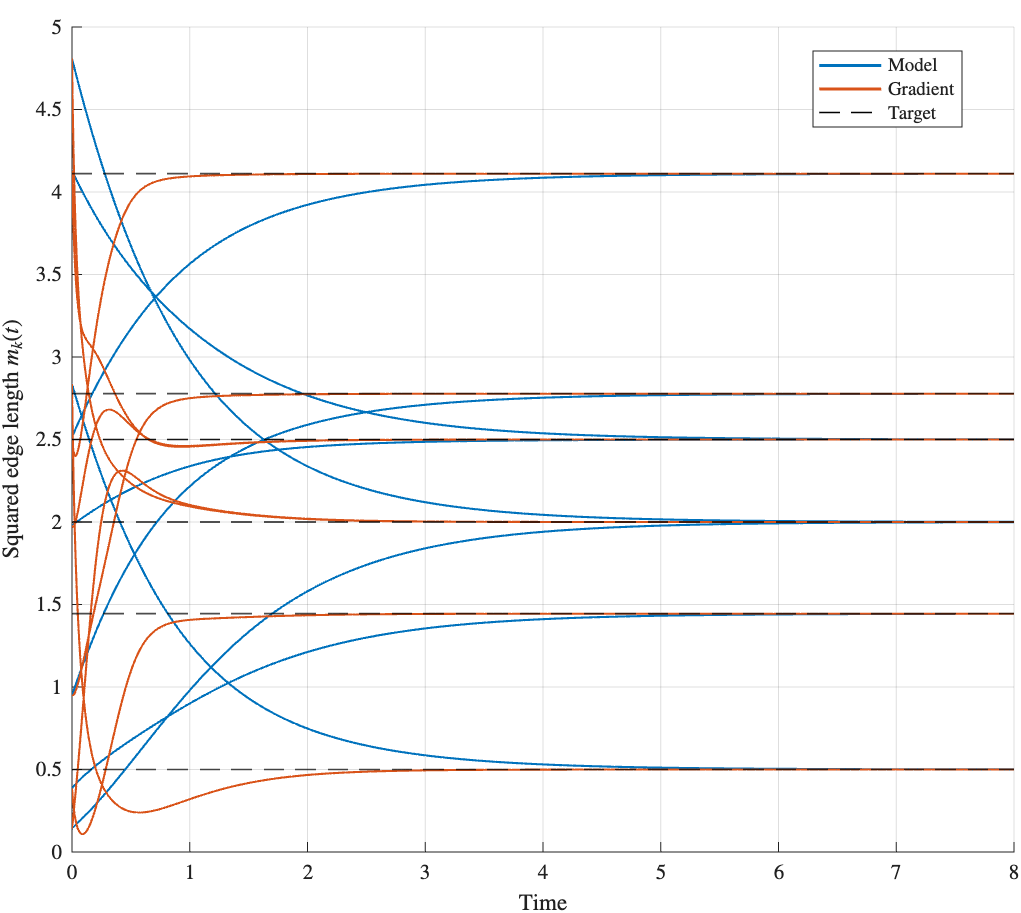}
    \caption{Squared edge trajectories under the model (blue) and
    gradient controllers (red).}
    \label{fig:edge_traj_compare}
  \end{subfigure}
  % (b) Node trajectories
  \begin{subfigure}[t]{0.24\textwidth}
    \centering
    \includegraphics[width=\linewidth]{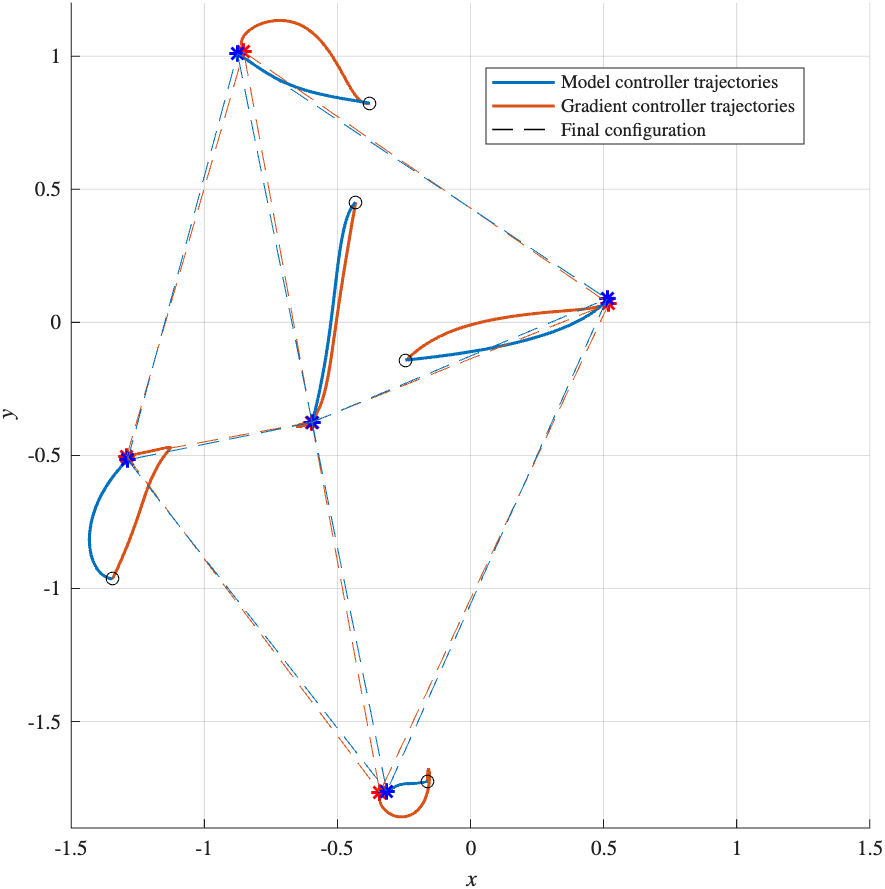}
    \caption{Node trajectories and final configurations for both
    controllers.}
    \label{fig:node_traj_compare}
  \end{subfigure}
  \caption{We compare
    the model controller
    \eqref{model_controller} and gradient controller \eqref{grad_Control}
  for $(W_5, p^*).$}
  \label{fig:model_vs_gradient}
\end{figure}

\begin{figure}[H]
  \centering
  \includegraphics[width=0.2\textwidth]{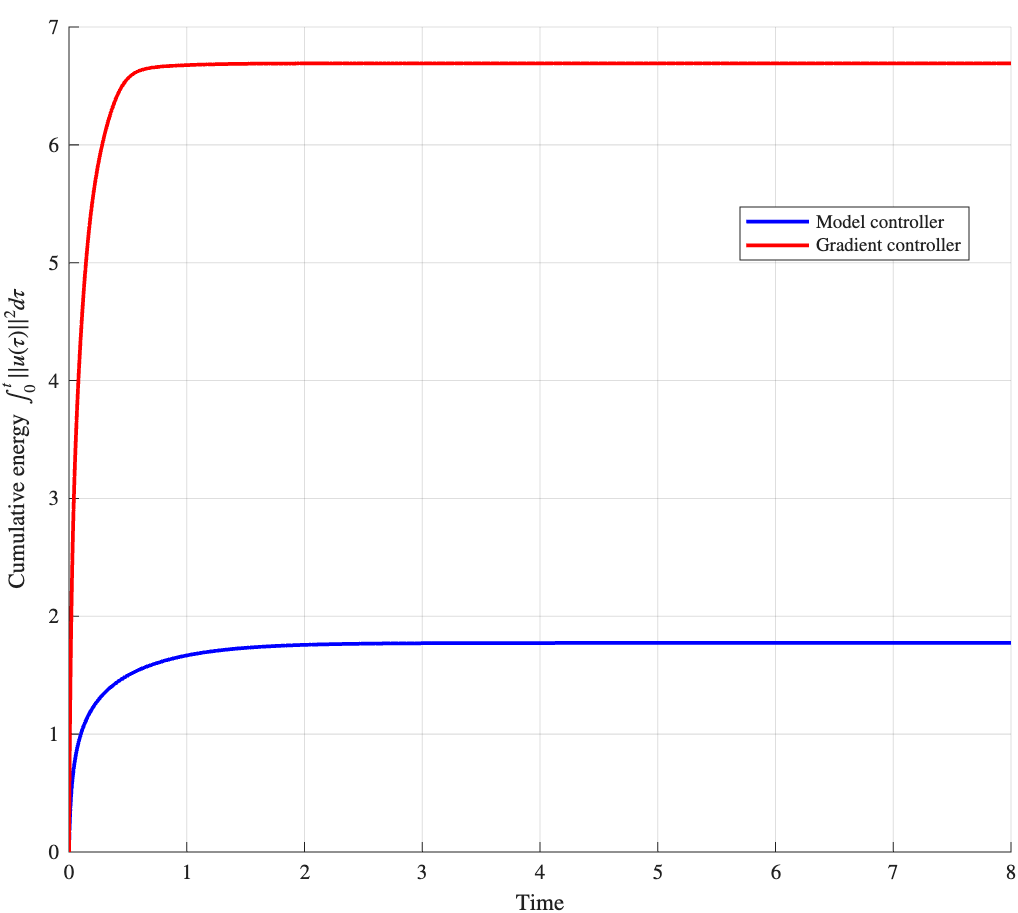}
  \caption{Cumulative control energy for the model controller (blue)
  and gradient controller (red) for $(W_5, p^*)$.}
  \label{fig:energy_traj_compare}
\end{figure}
We begin with the wheel graph $\mathrm W_5$ with node~5 as hub
(Figure~\ref{fig:wheel5_combined}(a)).
Figure~\ref{fig:model_vs_gradient} compares the model controller
\eqref{model_controller} with the gradient controller
\eqref{grad_Control}. Both achieve the desired distances
(Figure~\ref{fig:edge_traj_compare}), but the model controller
requires significantly less control energy
(Figure~\ref{fig:energy_traj_compare}), as predicted by
\eqref{eq:node-opt}. The resulting trajectories differ, but in
both cases the formation converges to a configuration congruent to
$p^*$ (Figure~\ref{fig:node_traj_compare}).

\subsubsection*{Target dependence of directed control}
\begin{figure}[H]
  \centering
  % (a) Edge trajectories
  \begin{subfigure}[t]{0.24\textwidth}
    \centering
    \includegraphics[width=\linewidth]{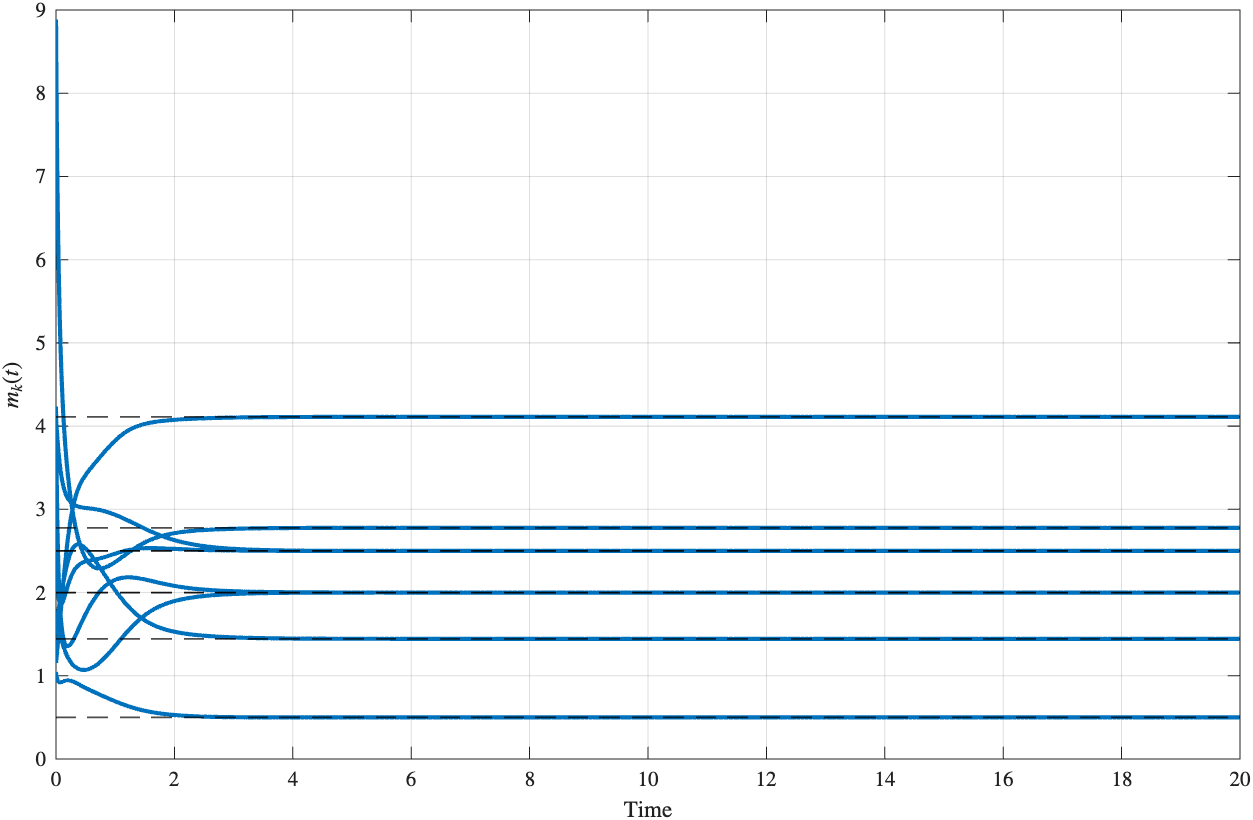}
    \caption{Squared edge trajectories under the directed controller.}
    \label{fig:dir_edge_traj}
  \end{subfigure}
  % (b) Node trajectories
  \begin{subfigure}[t]{0.24\textwidth}
    \centering
    \includegraphics[width=\linewidth]{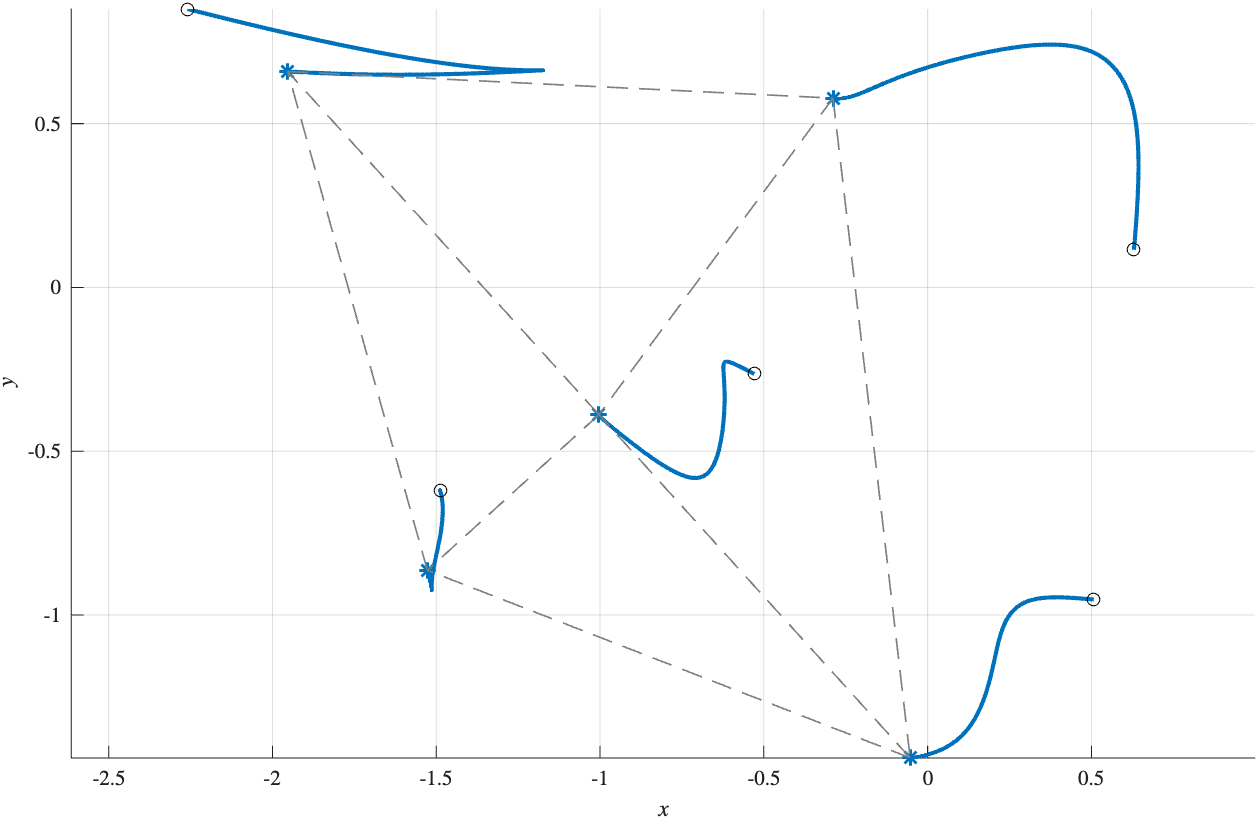}
    \caption{Node trajectories and final configuration.}
    \label{fig:dir_node_traj}
  \end{subfigure}
  \caption{Trajectories of the directed controller
    \eqref{dir_formation_ctr} for $({\small \overrightarrow{W_5}},p^*)$,
  which satisfies Theorem \ref{thm: dircontrol}.}
  \label{fig:good_dir_control}
\end{figure}

We next consider the directed wheel $\overrightarrow{\mathrm W}_5$
(Figure~\ref{fig:wheel5_combined}(b)). For target $p^*$, the directed
controller satisfies the sufficient condition in
Theorem~\ref{thm: dircontrol}, and the trajectories in
Figure~\ref{fig:good_dir_control} show exponential convergence.

\begin{figure}[H]
  \centering
  % (a) Edge trajectories
  \begin{subfigure}[t]{0.24\textwidth}
    \centering
    \includegraphics[width=\linewidth]{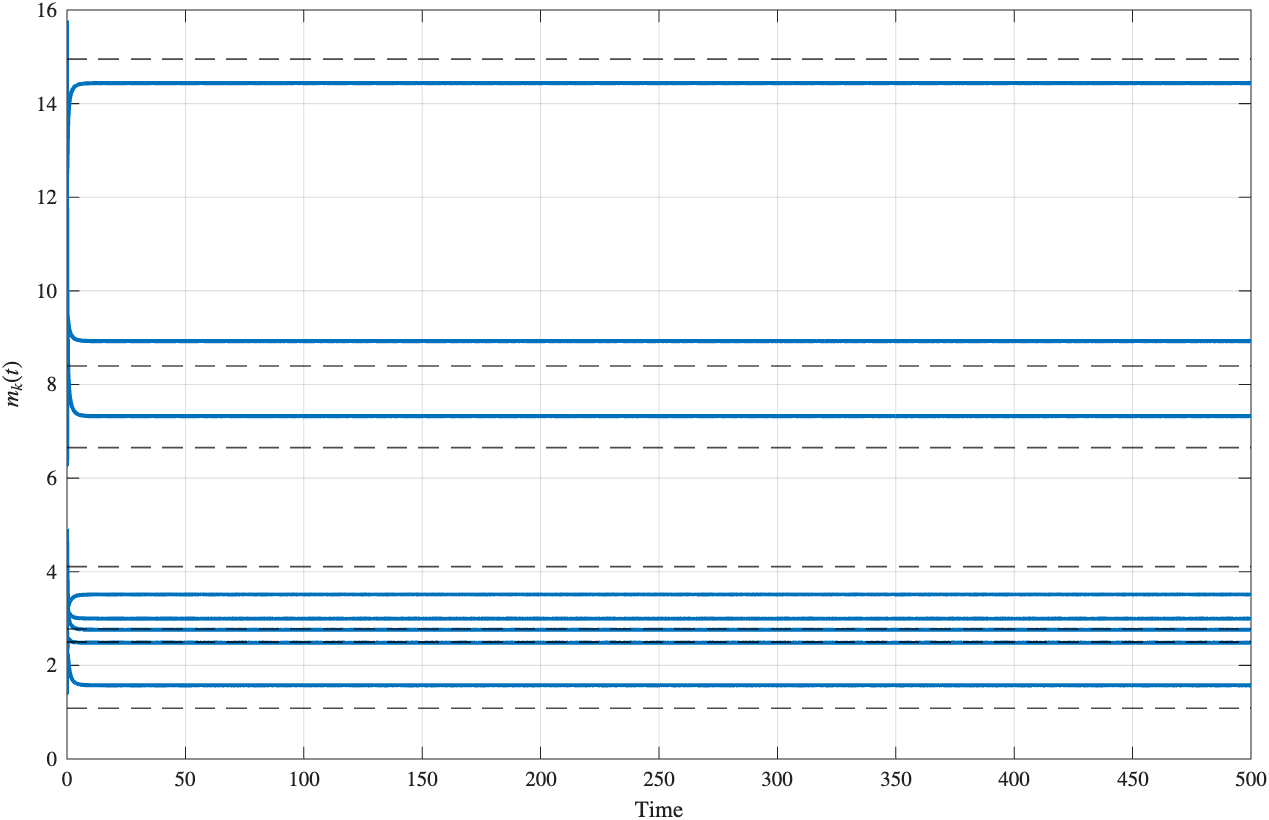}
    \caption{Squared edge trajectories under the directed controller.}
    \label{fig:baddir_edge_traj}
  \end{subfigure}
  % (b) Node trajectories
  \begin{subfigure}[t]{0.24\textwidth}
    \centering
    \includegraphics[width=\linewidth]{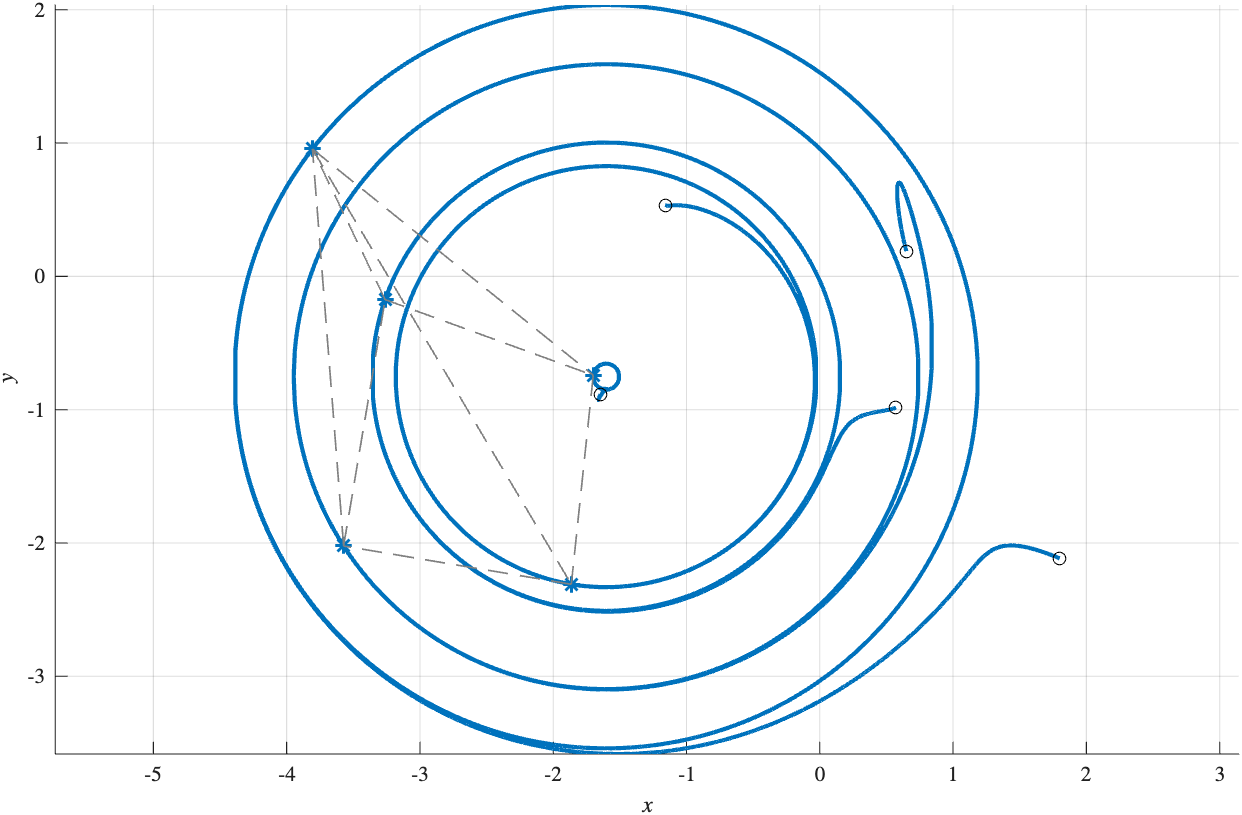}
    \caption{Node trajectories converge to a limit cycle.}
    \label{fig:baddir_node_traj}
  \end{subfigure}
  \caption{Trajectories of the directed controller
    \eqref{dir_formation_ctr} for $({\small \overrightarrow{W_5}},
    q^*),$ which does not
  satisfy Theorem \ref{thm: dircontrol}.}
  \label{fig:bad_dir_control}
\end{figure}

In contrast, with the \emph{same orientation} but different target $q^*$
(Figure~\ref{fig:wheel5_combined}(c)), the numerical test fails and
the system exhibits qualitatively different behavior: the edge
lengths converge, but the node trajectories approach a stable limit
cycle (Figure~\ref{fig:bad_dir_control}). Thus, stability depends
essentially on \emph{target geometry} and not merely on graph
orientation.

\subsubsection*{Persistence is not necessary}
\begin{figure}[H]
  \centering
  \begin{tikzpicture}

% Axis lines
    \draw[->, thin] (-1.9,0) -- (1.65,0) node[below right] {}; % x-axis
    \draw[->, thin] (0,-2.0) -- (0,1.75) node[above] {}; % y-axis

    \Vertex[size=.3, x=.11, y=-1.03, label=$1$, style={shape=circle}]{v1}
    \Vertex[size=.3, x=-.91,  y=-.11,  label=$2$, style={shape=circle}]{v2}
    \Vertex[size=.3, x=1.44,  y=1.64,  label=$3$, style={shape=circle}]{v3}
    \Vertex[size=.3, x=.35, y=-1.99,  label=$4$, style={shape=circle}]{v4}
    \Vertex[size=.3, x=-1.87,  y=1.53, label=$5$, style={shape=circle}]{v5}
    \Vertex[size=.3, x=1.61,  y=.77, label=$6$, style={shape=circle}]{v6}

    \Edge[Direct](v2)(v1)
    \Edge[Direct,bend=-15](v3)(v1)
    \Edge[Direct,bend=-15](v3)(v5)
    \Edge[Direct,bend=15](v4)(v2)
    \Edge[Direct,bend=5](v4)(v3)
    \Edge[Direct,bend=15](v5)(v1)
    \Edge[Direct](v5)(v6)
    \Edge[Direct](v6)(v2)
    \Edge[Direct,bend=15](v6)(v4)
    \Edge[Direct,bend=-15](v3)(v2)
    \Edge[Direct,bend=-15](v5)(v2)
\end{tikzpicture}
  \caption{A directed graph that is not persistent.}
  \label{fig:dirgraph}
\end{figure}
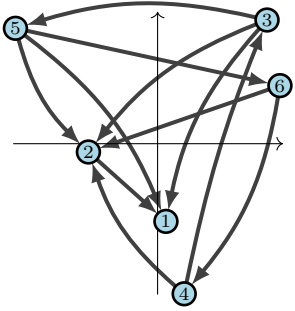
\begin{figure}[H]
  \centering
  % (a) Edge trajectories
  \begin{subfigure}[t]{0.24\textwidth}
    \centering
    \includegraphics[width=\linewidth]{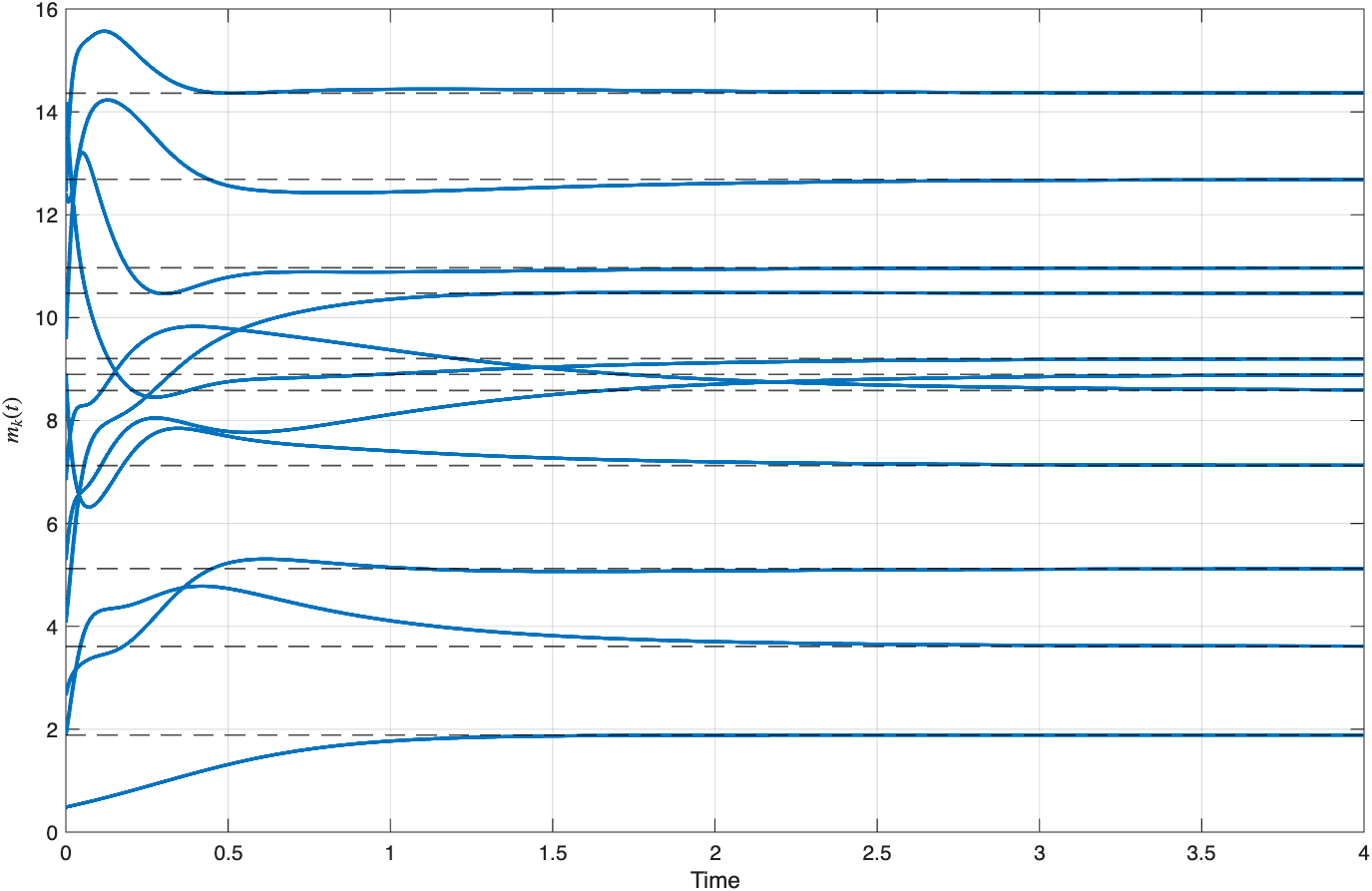}
    \caption{Squared edge trajectories under the directed controller.}
    \label{fig:nonpersistent_edge_traj}
  \end{subfigure}
  % (b) Node trajectories
  \begin{subfigure}[t]{0.24\textwidth}
    \centering
    \includegraphics[width=\linewidth]{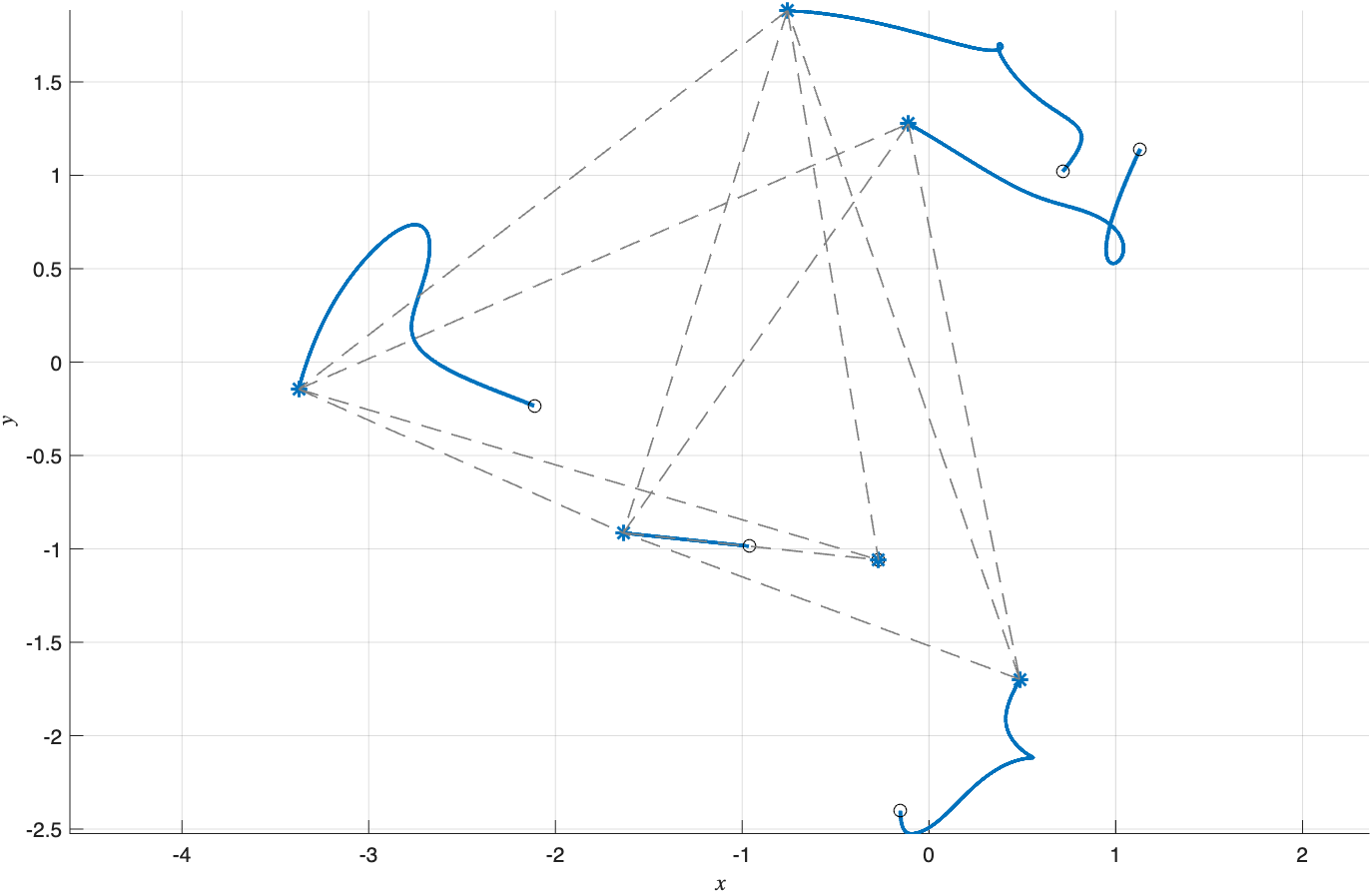}
    \caption{Node trajectories and the final configuration.}
    \label{fig:nonpersistent_node_traj}
  \end{subfigure}
  \caption{Trajectories of the directed controller
    \eqref{dir_formation_ctr} for $({\small
    \overrightarrow{\mathcal{G}}}, p^*)$, for a non-persistent graph
  ${\small \overrightarrow{\mathcal{G}}}$.}
  \label{fig:nonpersistent_control}
\end{figure}

Our final example uses the non-persistent directed graph
$\overrightarrow{\mathcal{G}}$ in Figure~\ref{fig:dirgraph}. Although
the directed graph fails the persistence criterion
\cite[Theorem~3]{Hendrickx_IJRNC2007}, the  underlying undericted
graph is $2$-rigid and the sufficient condition in
Theorem~\ref{thm: dircontrol} is satisfied. As shown in
Figure~\ref{fig:nonpersistent_control}, the directed controller again
converges.

\section{After persistence: admissibility}\label{sec.combinatorics}
We have shown that persistence does not describe the behavior of
the directed controller \eqref{dir_formation_ctr}.
Indeed, Figure \ref{fig:bad_dir_control} shows that it is not
sufficient (not
sensitive to the target geometry), while Figure
\ref{fig:nonpersistent_control}
shows that persistence is not necessary (not linked to the geometry and
dynamics).
We now describe an alternative inspired by the algebraic shadow
of the dynamics in the family \eqref{eq:pm-dynamics}.

We begin with the linearized edge dynamics of a controller in the family
\eqref{eq:pm-dynamics} at a target configuration
$(p^*,m^*)\in \mathcal{S}$.  Considering a perturbation
$(\delta p,\delta m)\in T_{(p^*,m^*)} \mathcal{S}$ and using the
notation $e = m^* - m$ for the edge vector, we get $\delta e =
-\delta m$, so
$\dot{\delta e} = -\eta\delta m$.  Hence, the linearized edge dynamics
on $T_{m^*} \mathcal{Q}$ are described by the linear ODE
\begin{equation}\label{eq: linear dynamics}
  \dot e = -\eta^*|_{T_{m^*} \mathcal{Q}}e.
\end{equation}
From the linear dynamics at the target we define two
generic properties of such a controller.
\begin{definition}\label{def: admissible}
  A controller  \eqref{eq:pm-dynamics} defined
  by rational functions is called \emph{dynamically admissible} if,
  for every generic
  target configuration $(p^*,m^*)\in \mathcal{S}$, $\eta^*|_{T_{m^*}
  \mathcal{Q}}$
  is hyperbolic.  It is called \emph{algebraically admissible} if
  $\eta^*|_{T_{m^*} \mathcal{Q}}$ is invertible for generic targets.
\end{definition}
Dynamic admissibility (all eigenvalues have non-zero real parts)
implies algebraic admissibility (no eigenvalue is zero).
Admissibility is our proposed replacement
for persistence when studying the directed controller
\eqref{dir_formation_ctr}. Both types
of admissibility are easy to check with a randomized algorithm: pick a
random configuration $p^*$ and then check the eigenvalues of
$\eta_{(p^*,F(p^*))}$. Admissibility can be checked
independently of the target because it concerns generic eigenstructure,
as opposed to the certificates in Theorems \ref{thm: main} and
\ref{thm: dircontrol}, which are about Lyapunov geometry near the
target.  We also note that the test for both types of admissibility
is polynomial time, while the test for persistence described
by \cite{Hendrickx_IJRNC2007} is exponential in the number
of redundant edges, and, when $d\ge 3$, also uses linear algebra.

We finish by connecting admissibility to the convergence properties
of our controllers.
\begin{theorem}[Admissibility is necessary]\label{thm: admiss}
  Let \\Assumption~\ref{asump.drigid-geneic} hold.
  If a controller defined by
  \eqref{eq:pm-dynamics-a}--\eqref{eq:pm-dynamics-b} defined
  by rational functions is locally exponentially stable for some
  target formation $(p^*,m^*)\in \mathcal{S}$, then it is dynamically
  admissible.
\end{theorem}
\begin{pf}
  By Corollary 4.3 in \cite{khalil} (used in an appropriate chart on
  $\mathcal{Q}$),
  if a non-linear IVP converges exponentially, the linearization of the
  dynamics at the equilibrium must be Hurwitz, and, hence
  hyperbolic.  From \eqref{eq: linear dynamics}, this implies that
  $\eta^*|_{T_{m^*}\mathcal{Q}}$ is hyperbolic.  Because dynamic
  admissbility is a
  generic property, one example where it is satisfied implies the
  result.\hfill$\qed$
\end{pf}
In particular, dynamic and algebraic admissibility are \emph{necessary} for
local exponential convergence of any control law in our family.  This is
in parallel to the one direction of the role of generic $d$-rigidity for
the model and undirected controllers: a combinatorial constraint that any
locally stable design must satisfy.

\section{Further directions}
Our work suggests some immediate open problems.  A combinatorial
classification of the
dynamically and algebraically admissible orientations of a graph is a
natural next step,
and it would be interesting to identify special families of
$d$-rigid graphs and
orientations for which the sufficient condition in Theorem
\ref{thm: dircontrol}
holds generically.  It would also be of interest to know whether
the sufficient
stability condition in Theorem \ref{thm: dircontrol} is also necessary,
as it is in Theorem \ref{thm: model}.

\section*{Acknowledgements} We thank Ashley Bao, Martin Deraas,
Shlomo Gortler, Karime Hernandez, Tony Nixon, Bernd Schulze and
Audrey St John for
helpful conversations and feedback.

\end{document}